\documentclass[12pt]{article}
\usepackage{isolatin1,times,a4wide,doublespace,amsmath,latexsym,amssymb}
\newtheorem{Thm}{Theorem}[section]

\newtheorem{Cor}[Thm]{Corollary}
\newtheorem{Lem}[Thm]{Lemma}
\newtheorem{Prop}[Thm]{Proposition}
\newtheorem{Def}[Thm]{Definition}

\newtheorem{Not}[Thm]{Notation}

\newcommand{\bN}{{\mathbb N}}
\newcommand{\bR}{{\mathbb R}}

\newcommand{\bproof}{\noindent{\bf Proof: }}
\newcommand{\eproof}{\hfill $\Box$\\}
\newcommand{\bremark}{\noindent{\bf Remark: }}
\newcommand{\eremark}{\hfill \\}

\newcommand{\cS}{{\mathcal S}}
\def\ldots{\mathinner{\ldotp\ldotp\ldotp}}
\def\ldots{\mathinner{\cdotp\cdotp\cdotp}}

\begin{document}

\title{A Banach space with a symmetric basis which is of weak cotype 2
but not of cotype 2}
\author{Peter G. Casazza\footnote{Supported by NSF grant DMS 9706108.}
\and Niels J. Nielsen\footnote{Supported by the Danish Natural Science
Research
Council, grant 9801867.}}
\date{}
\maketitle

\begin{abstract}
We prove that the symmetric convexified Tsirelson space
is of weak cotype 2 but not of cotype 2.
\end{abstract}

\section*{Introduction}

Weak type 2 and weak cotype 2 spaces were originally introduced  and
investigated by V.D.\ Milman and G.\ Pisier in \cite{MP} and weak
Hilbert
spaces by Pisier in \cite{P1}. A
further detailed investigation can be found in Pisier's book
\cite{P}. The first example of a weak Hilbert space which is not
isomorphic to a Hilbert space is the 2-convexified Tsirelson space
(called the convexified Tsirelson space in this paper). This follows
from
the results of W.B.\ Johnson in \cite{J}. For a detailed study of the
original Tsirelson space we refer to \cite{CS}.

Let $X$ be a Banach space with a symmetric basis. It was proved in
\cite{P} that if $X$ is a weak Hilbert space, then it is isomorphic to a
Hilbert space and this has lead to the belief that if $X$ is just of
weak cotype
2, then it is of cotype 2. However, this turns out not necessarily to be
the
case. The main result of this paper states that the symmetric
convexified
Tsirelson space is of weak cotype 2  but not of cotype 2.

We now wish to discuss the arrangement of this paper in greater detail.

In Section 1 we give some basic facts on properties related to weak type
2 and weak cotype 2 while Section 2 is devoted to a review of some
results on the convexified Tsirelson space which we need for our main
result. Most of these results are stated without proofs since they can
be proved in a similar manner as the corresponding results for the
original Tsirelson space.

In Section 3 we make the construction of the symmetric convexified
Tsirelson space, investigate its basic properties and prove our main
result stated above.

\section*{Acknowledgement}

The authors are indebted to Nigel Kalton for communicating Theorem
\ref{Nigel} to us.

\section{Notation and Preliminaries}\label{PR}
\setcounter{equation}{0}

In this paper we shall use the notation and terminology commonly used in
Banach space theory as it appears in \cite{LT1}, \cite{LT2} and
\cite{T}. $B_X$ shall always denote the closed unit ball of the Banach
space $X$ and if $X$ and $Y$ are Banach spaces, then $B(X,Y)$
($B(X)=B(X,X)$) denotes the
space of all bounded linear operators from $X$ to $Y$.

We let $(g_n)$ denote a sequence of independent standard Gaussian
variables on a fixed probability space $(\Omega,\cS,\mu)$ and recall
that a Banach space $X$ is said 
to be of type p, $1\le p \le 2$, (respectively cotype p, $2\le p <
\infty$) if there is a constant $K\ge 1$ so that for all finite sets
$\{x_1,x_2,\dots,x_n\}\subseteq X$ we have
\begin{equation}
\label{eq1.1}   
\big(\int \big\|\sum_{j=1}^n g_j(t)x_j\big\|^p
d\mu(t)\big)^\frac1p \le K \big(\sum_{j=1}^n \|x_j\|^p \big)^\frac1p
\end{equation}
(respectively
\begin{equation}
\label{eq1.2}
K \big(\sum_{j=1}^n \|x_j\|^p \big)^\frac1p \le\big(\int
\big\|\sum_{j=1}^n g_j
(t)x_j\big\|^pd\mu(t)\big)^\frac1p).
\end{equation}

The smallest constant $K$ which can be used in (\ref{eq1.1})
(respectively (\ref{eq1.2})) is denoted by $K^p(X)$ (respectively
$K_p(X)$).

If $L$ is a Banach lattice and $1\le p< \infty$, then $L$
is said to be $p$-convex (respectively $p$-concave) if there
is a constant $C\ge 1$ so that for all finite sets
$\{x_1,x_2,\dots,x_n\}\subseteq L$ we have
\begin{equation}
\label{eq1.3}
\|(\sum^n_{j=1} |x_j|^p)^{\frac1p}\| \le C(\sum^n_{j=1}
\|x_j\|^p)^{\frac1p}
\end{equation}
(respectively
\begin{equation}
\label{eq1.4}
(\sum^n_{j=1} \|x_j\|^p)^{\frac1p} \le C\|(\sum^n_{j=1}
|x_j|^p)^{\frac1p}\| ).
\end{equation}

The smallest constant $C$ which can be used in (\ref{eq1.3})
(respectively (\ref{eq1.4})) is denoted by $C^p(L)$ (respectively
$C_p(L)$). 

It follows from \cite[1.d.6 (i)]{LT2} that if $L$ is of finite concavity
(equivalently of finite cotype), then there is a constant $K \ge 1$ so
that
\begin{equation}
\label{eq1.5}
\frac{1}{K}\|(\sum^n_{j=1} |x_j|^2)^{\frac12}\| \le \big(\int
\big\|\sum_{j=1}^n g_j(t)x_j\big\|^2
d\mu(t)\big)^\frac12 \le K\|(\sum^n_{j=1} |x_j|^2)^{\frac12}\|
\end{equation} 

A Banach space $X$ is said to be of weak type 2 if there is a constant
$C$ and a $\delta$, $0<\delta<1$, so that whenever $E\subseteq X$ is a
subspace, $n\in\bN$ and $T\in B(E,\ell_2^n)$, then there is an
orthogonal
projection $P$ on $\ell_2^n$ of rank larger than $\delta n$ and an
operator $S\in B(X,\ell_2^n)$ with $Sx=PTx$ for all $x\in E$ and
$\|S\|\le C\|T\|$.

Similarly $X$ is called a weak cotype 2 if there is a constant $C$ and a
$\delta$, $0<\delta<1$, so that whenever $E\subseteq X$ is a finite
dimensional subspace, then there is a subspace $F\subseteq E$ so that
$\dim F\ge \delta \dim E$ and $d(F,\ell_2^{\dim F})\le C$.

Our definitions of weak type 2 and weak cotype 2 space are not the
original ones, but are chosen out of the many equivalent
characterizations given by Pisier \cite{P}.

A weak Hilbert space is a space which is both of weak type 2 and weak
cotype 2.

If $A$ is a set we let $|A|$ denote the cardinality of $A$.
\begin{Def}
If $(x_{n})$ and $(y_{n})$ are sequences in a Banach space
$X$, we say that $(x_{n})$ is dominated by $(y_{n})$ if
there is a constant $K>0$ so that for all finitely non-zero
sequences of scalars $(a_{n})$ we have
$$
\|\sum_{n}a_{n}x_{n}\|\le K\|\sum_{n}a_{n}y_{n}\|.
$$
\end{Def}

We will need some information about property $(H)$ and related
properties.

\begin{Def}\label{D1.2}
A Banach space $X$ has property $(H_{2})$ if there is a function
$C(\cdot , \cdot )$ so that 
for every $0<\delta
<1$ and for every normalized $\lambda$-unconditional basic sequence
$(x_{i})_{i=1}^{n}$ in $X$ there is a subset $F\subseteq\bN$ 
such that $|F|\ge {\delta}n$ and $(x_{i})_{i\in F}$ is 
$C({\lambda},{\delta})$-equivalent to the unit vectors basis of
${\ell}_{2}^{|F|}$.  If we only have that $(x_{i})_{i\in F}$ is
$C({\lambda},{\delta})$ dominated by the unit vector basis of
${\ell}_{2}^{|F|}$, we say that $X$ has property upper
$(H_{2})$.  Similarly, we define property lower $(H_{2})$.
\end{Def}

\begin{Def}\label{D1.3}
A Banach space $X$ is said to have property $(H)$ if there is
a function $f(\cdot )$ so that for every 
normalized $\lambda$-unconditional
basic sequence $(x_{i})_{i=1}^{n}$ in $X$, we have
$$
\frac{1}{f(\lambda )}n^{1/2} \le \|\sum_{i=1}^{n}x_{i}\| \le
f(\lambda )n^{1/2}.
$$
Similarly, we can define property upper $(H)$ and property
lower $(H)$.
\end{Def}

The following is clear.

\begin{Prop}\label{P20}
Property upper (resp. lower) $(H_{2})$ implies upper (resp.
lower) $(H)$.
\end{Prop}

We will see later that the converses of Proposition \ref{P20}
fail.

The next result shows that any percentage of the basis will
work in the definition of $(H_{2})$.  The proof follows from
the argument of Pisier
\cite[Proposition 12.4, page 193]{P}.    

\begin{Lem}\label{L1.4}
For a Banach space $X$, the following are equivalent:

(1)  $X$ has property upper (resp. lower) $(H_{2})$.

(2)  There exists one $0<\delta <1$ satisfying the conclusion
of property upper (resp. lower) $(H_{2})$.
\end{Lem}

The corresponding result for property $(H)$ is in
\cite[Proposition Ae1, page 14]{CS}.

\begin{Lem}\label{L1.5}
For a Banach space $X$, the following are equivalent:

(1)  $X$ has property upper (resp. lower) $(H)$.

(2)  There is a $0<{\delta}<1$ so that for every 
$\lambda$-unconditional basic sequence $(x_{i})_{i=1}^{n}$
in $X$ there is a subset $F\subset \{1,2,\ldots , n\}$
with $|F|\ge {\delta}n$ and $(x_{i})_{i\in F}$ has property
upper (resp. lower) $(H)$.
\end{Lem}

The next theorem is due to Pisier \cite[Proposition 12.4]{P}.

\begin{Prop}\label{P1.6}
Every weak Hilbert space has property $(H_{2})$.
\end{Prop}

We also have from Pisier \cite[Proposition 10.8, page 160 and
Proposition 11.9, page 174]{P}:

\begin{Prop}\label{P21}
The following implications hold for a Banach space $X$:

(1)  Weak cotype 2 implies property lower $(H)$.

(2)  Weak type 2 implies property upper $(H)$.
\end{Prop}

The converses of Proposition \ref{P21} are open questions.  However,
for Banach lattices it is known that property $(H)$, property $(H_{2})$
and being a weak Hilbert space are all equivalent.  This is a
result of Nielsen and Tomczak-Jaegermann \cite{NT}.

\section{Convexified Tsirelson Space}\label{CT}
\setcounter{equation}{0}

Since there is only a ``partial theory'' developed for 
the convexified Tsirelson space $T^2$, we will review what we
need here.

\begin{Not}
If $E,F$ are sets of natural numbers, we write $E<F$ if
for every $n\in E$ and every $m\in F$, $n<m$.  If $E=\{k\}$,
we just write $k<F$ for $E<F$.
\end{Not} 

\begin{Def}
We define the convexified Tsirelson space $T^{2}$ as the set of
vectors $x=\sum_{n}a_{n}t_{n}$ for which the recursively
defined norm below is finite.  
\begin{equation}
\label{eq2.1}
\|x\|_{T^{2}} = \mbox{max}\{\mbox{sup} |a_{n}|,
2^{-1/2}\mbox{sup}\left ( \sum_{j=1}^{k}\|E_{j}x\|^{2} \right ) ^{1/2}
\},
\end{equation}
where the second ``sup'' is taken over all choices 
$$
k\le E_{1} < E_{2} < \cdots < E_{k},
$$
and $Ex = \sum_{n\in E}a_{n}t_{n}$.
\end{Def}

We will now list the known results for this space 
(which we will need)  and where
they can be found.  The first result can be found in
\cite{CS} and \cite{P}.

\begin{Prop}\label{P2.2}
The unit vectors $(t_{n})$ form a 1-unconditional basis for
$T^2$.  The space $T^2$ is of type 2 and weak cotype 2 but
does not contain a Hilbert space.
\end{Prop}

Next we need to see which subsequences of the unit vector
basis of $T^{2}$ are equivalent to the original basis. To do this we
need:

\begin{Not}\label{N2.7}
The fast growing hierarchy from logic is a family
of functions on $\mathbb N$ given by:  $g_{o}(n) =
n+1$, and for $i\ge 0$, $g_{i+1}(n) = g_{i}^{(n)}(n)$, where
for any function $f$, $f^{(n)}$ is the n-fold iteration of $f$.
We also set $\text{exp}_{0}(n) = n$ and for $i,n\ge 1$,
$$
\text{exp}_{i}(n) = 2^{\text{exp}_{i-1}(n)}.
$$
Finally we let $\text{log}_{0}(n) = n$, and for n large
enough so that $\text{log}_{i-1}(n)>0$, let
$$
\text{log}_{i}(n) = \text{log(log}_{i-1}(n)).
$$ 
\end{Not}

The next result is due to Bellenot \cite{B}.  He does this
result in the original Tsirelson's space $T$, but the proof
works perfectly well in $T^{2}$.

\begin{Prop}\label{P2.8}
A subsequence $(t_{k_{n}})$ of $(t_{n})$ is equivalent to
$(t_{n})$ if and only if there is a natural number $i$ so
that $k_{n}\le g_{i}(n)$, for all large n.  Moreover,
$(t_{k_{n}})$ always 1-dominates $(t_{n})$ and
there is a constant $K\ge 1$ so that the equivalence
constant is $K^{i}$ for the case $g_{i}(n)$.
\end{Prop}

One important consequence is (see Pisier \cite{P} or
Casazza and Shura \cite{CS}). 

\begin{Prop}\label{P2.11}
Every $g_{i}(n)$-dimensional subspace of span $(t_{j})_{j\ge n}$
is $K^{i}$-isomorphic to a Hilbert space and $K^{i}$-complemented
in $T^{2}$.
\end{Prop}

If $X$ is a weak Hilbert space with an unconditional basis, then it
follows from \cite{NT} that the conclusion of Proposition \ref{P2.11}
remains true after a suitable permutation of the basis.

The next result comes from \cite[Theorem IV.b.3, page 39]{CS}.
The theorem there is proved for the regular Tsirelson space but
the techniques easily adapt to convexified space.

\begin{Prop}\label{P2.10}
Every n-dimensional subspace of $T^{2}$ is $K^{i}\text{log}_{i}(n)$
isomorphic to ${\ell}_{2}^{n}$.
\end{Prop}

We need one more result on convexified Tsirelson.

\begin{Prop}\label{P2.12}
If $x = \sum_{j}a_{j}t_{j}\in T^2$, then for all $n\in \mathbb N$,
$$
\|\sum_{j}a_{j}t_{nj}\|_{T^2}\le 2K^{i}({\text log}_{i}n)\|x\|_{T^2}.
$$
\end{Prop}

\bproof
By Proposition \ref{P2.8} and Proposition \ref{P2.10} we have

\begin{eqnarray*}
\|\sum_{j}a_{j}t_{nj}\|_{T^2} \le \|\sum_{j=1}^{n}a_{j}t_{nj}
\|_{T^2} + \|\sum_{j=n+1}^{\infty}a_{j}t_{nj}\|_{T^2} 
& \le & \left ( \sum_{j=1}^{n}|a_{j}|^{2} \right )^{1/2} +
\|\sum_{j=n+1}^{\infty}a_{j}t_{j}\|_{T^2} \\
\le K^{i}({\text log}_{i}n)\|\sum_{j=1}^{n}a_{j}t_{j}\| + K
\|\sum_{j=n+1}^{\infty}a_{j}t_{j}\|_{T^2} 
&\le & 2K^{i}({\text log}_{i}n)\|x\|.
\end{eqnarray*}
\eproof 

\section{Symmetric Convexified Tsirelson Space}\label{SCT}
\setcounter{equation}{0}

There is almost no existing theory for the symmetric convexified
Tsirelson space.  But there is a theory for the symmetric Tsirelson
space.  We will list the results we need on this topic. They can be
found in
Casazza and Shura \cite[ Chapter X.E]{CS}.   

\begin{Not}\label{N3.1}
For $T^{2}$ or $(T^{2})^{*}$ we will work with the non-decreasing
rearrangement operator $D$.  That is, if $x = \sum_{n}a_{n}t_{n}$
then $Dx = \sum_{n}a_{n}^{*}t_{n}$ where $(a_{n}^{*})$ is the
non-decreasing re-arrangement of the non-zero $a_{n}'s$ where
by non-decreasing we mean the absolute values are non-decreasing.  
\end{Not}

The construction of Chapter VIII of \cite[Chapters VIII and X.B]{CS}
shows

\begin{Prop}\label{P3.2}
Let $\Pi$ denote the group of all permutations of $\mathbb N$.
There is a constant $K\ge 1$ so that for any $x=\sum_{n}
a_{n}t_{n}^{*}
\in (T^{2})^{*}$ we have
\begin{equation}
\label{eq3.1}
\|x\|_{s^{*}} =:
{\text sup}_{\sigma \in \Pi}\|\sum_{n}a_{{\sigma}(n)}
t_{n}^{*}\| \le K \|Dx\| \le 
{K}{\text sup}_{\sigma \in \Pi}\|\sum_{n}a_{{\sigma}(n)}
t_{n}^{*}\|.
\end{equation}

\end{Prop}

We will define the {\em dual space of the symmetric convexified
Tsirelson
space} first because it is natural in terms of the above.

\begin{Def}\label{D3.3}
We let $S[(T^{2})^{*}]$ be the family of all vectors for which
$\|x\|_{s^{*}}$ is finite.  Then this is a Banach space with
a natural symmetric basis, denoted $(t_{n}^{s*})$, called the
dual space of the symmetric convexified Tsirelson space.
\end{Def}

To define the {\em the symmetric convexified Tsirelson space} we need a
result
kindly communicated to us by N.J. Kalton.

Let $X$ be a Banach sequence space.  Define the permutation operators
$S_{\sigma}(\xi)= (\xi_{\sigma(n)})_{n=1}^{\infty}$ for $\sigma\in \Pi$
and let $L_k^j$ to be the linear map such that $L_k^j(e_n)=e_{kn+j}$ for
all $n\in \bN$. Finally we let $c_{00}$ denote the spaces af real
sequences which are eventually 0.

\begin{Thm} \label{Nigel}
Suppose $X$ is a Banach sequence space which is $p$-convex
and $q$-concave where $1<p<q<\infty.$  Suppose
 $\max_{0\le j<k}\|L_k^j\|\le Ck^a$ where
$a+p^{-1}<1.$  Then
$$ \|\xi\|_{X_{inf}}=\inf_{\sigma\in \Pi}\|S_{\sigma}\xi\|_X, \ x\in
c_{00}$$ defines
a quasi-norm on $c_{00}$ which is equivalent to a norm.  The dual of $
X_{inf}$ is $X^*_{sup}$ where
$$ \|\xi\|_{X^*_{sup}}= \sup_{\sigma\in\Pi} \|S_{\sigma}\xi\|_{X^*}.$$
\end{Thm}

\bproof
Let us start by supposing $x_1,\ldots,x_k\in c_{00}$ are disjointly
supported and that $\sigma_1,\ldots,\sigma_k\in\Pi.$
Then
\begin{align*} \|x_1+\cdots+x_k\|_{X_{inf}} &\le
\|\sum_{j=1}^kL^{j-1}_kS_{\sigma_j}x_j\|_X\\
&\le (\sum_{j=1}^k\|L^{j-1}_kS_{\sigma_j}x_j\|_X^p)^{\frac1p}\\
&\le Ck^a (\sum_{j=1}^k\|S_{\sigma_j}x_j\|_X^p)^{\frac1p}.\end{align*}
Now taking an infimum over $\sigma_j$ gives
\begin{equation}\label{disjoint} \|x_1+\cdots+x_k\|_{X_{inf}} \le Ck^a
(\sum_{j=1}^k\|x_j\|_{X_{\inf}}^p)^{\frac1p}.\end{equation}

Let us use (\ref{disjoint}) first to show that $\|\cdot\|_{X_{inf}}$ is
a quasi-norm.  Indeed if $x,y\in c_{00}$ then
$$ \|x+y\|_{X_{inf}}\le 2\|\max(|x|,|y|)\|_{X_{inf}}\le
2^{a+1}C(\|x\|_{X_{inf}}+\|y\|_{X_{inf}}).$$

Next note that (\ref{disjoint}) implies
$$ \|x_1+\cdots+x_k\|_{X_{inf}} \le Ck^{a+\frac1p}\max_{1\le j\le
k}\|x_j\|_{X_{inf}}.$$  From this it follows easily that if
$a+\frac1p<\frac1r<1$ we have
$$ \|x_1+\cdots+x_k\|_{X_{inf}} \le
C_r(\sum_{j=1}^k\|x_j\|^r)^{\frac1r}$$ for disjoint $x_1,\ldots,x_k.$
Thus we have an upper $r$-estimate for $X_{inf}.$

It is trivial to show $X_{inf}$ has a lower $q$-estimate.  Now
by \cite[Theorem 4.1]{K1} (a simpler proof is given in \cite[Theorem
3.2]{K2}
) it follows that
$X_{inf}$ is {\em lattice-convex} and this means that an upper
$r$-estimate implies (lattice) $s$-convexity for all $s<r$ (Theorem 2.2
of \cite{K1}). Hence
$X_{inf}$ is $r$-convex for every $r$ with $a+\frac1p<\frac1r.$  In
particular $1$-convexity implies the quasi-norm is equivalent to a norm.
In fact $X^*_{inf}$ is a reflexive Banach space.

Now it is obvious that $X_{inf}\subset (X^*_{sup})^*$ and
$X^*_{sup}\subset
(X_{inf})^*.$  Hence it follows easily that
$(X_{inf})^*=X^*_{sup}.$
\eproof

\bremark
We can apply the above result to the case of the weighted
$\ell_p-$space $X$, with $1<p<\infty$ defined by the norm
$$ \|\xi\|_X=(\sum_{n=1}^{\infty}|\xi_n|^pw_n)^{\frac1p}$$ where
$(w_n)$ is an increasing sequence satisfying an estimate of the form
$$ w_{kn}\le Ck^aw_n$$ where $a<p-1.$  The $X_{inf}$ is defined by
the quasi-norm
$$ \|\xi\|_{X_{inf}}=(\sum_{n=1}^{\infty}(\xi_n^*)^pw_n^p)^{\frac1p}$$
where $(\xi_n^*)$ is the decreasing rearrangement of $(|\xi_n|).$
In this case $X_{sup}$ is the Lorentz space $d((w_n)^{-q/p},q).$

This result can be rephrased.  If $(v_n)$ is a positive decreasing
sequence satisfying an estimate $v_n\le Ck^bv_{kn}$ where $b<1$ then
$d((v_n),p)^* $ can be identified with the space of all sequences
$(\xi_n)$ so that
 $$ (\sum_{n=1}^{\infty}(\xi^*_n)^qv_n^{-q/p})^{\frac1q}<\infty.$$
This result is a special case of results of Reisner \cite{R}.
\eremark

Proposition VIII.a.8 of \cite{CS} states that the decreasing
rearrangement
operator $D$ is a bounded non-linear operator on the original Tsirelson
space
$T$.
This result then immediately carries over to the 2-convexification of
$T$ which is our convexified Tsirelson space $T^2$.  By Proposition
\ref{P2.12} we have that Theorem \ref{Nigel} holds in this case.  
We summarize this in the following result:

\begin{Prop}\label{P3.4}
There is a constant $K\ge 1$ so that for any $x= \sum_{n}a_{n}
t_{n} \in T^{2}$ we have
\begin{equation}
\label{eq3.2}
 {\text inf}_{\sigma \in \Pi}
\|\sum_{n}a_{{\sigma}(n)}
t_{n}\| \le \|Dx\| \le 
{K}{\text inf}_{\sigma \in \Pi}\|\sum_{n}a_{{\sigma}(n)}
t_{n}\| .
\end{equation}

Moreover, there is a norm $\|\cdot \|_{s}$
on the set of vectors for which $\|Dx\|< \infty$ satisfying
\begin{equation}
\label{eq3.3}
\frac{1}{K}\|x\|_{s}\le \|Dx\|\le K\|x\|_{s}.
\end{equation}
\end{Prop}

Note that our operator $D$ does not satisfy
a triangle inequality, but does with the constant $K$ on
the sum side of the triangle inequality.

\begin{Def}\label{D3.5}
The symmetric convexified Tsirelson space is the Banach space
$S(T^{2})$ of vectors for which $\|x\|_{s}<\infty$ with natural
unit vector basis $(t_{n}^{s})$.  
By Theorem \ref{Nigel} this is a reflexive Banach space
whose dual space is $S[(T^{2})^{*}]$.   
\end{Def}

It is known \cite{CS} that every infinite dimensional subspace
of $S(T^{2}))$ contains a subspace which embeds into $T^{2}$.
In particular $S(T^{2})$ is a Banach space with a natural
symmetric basis which has no subspaces isomorphic to $c_{0}$ or
${\ell}_{p}$ for $1\le p<\infty$.  Also $T^{2}$ embeds into
$S(T^{2})$.  
Since the unit vector basis of ${\ell}_{2}$ uniformly dominates
all block bases of $(t_{n})$ in $T^{2}$, it follows that the
unit vector basis of $S(T^{2})$ is also dominated by the 
unit vector basis of ${\ell}_{2}$.

\begin{Prop}\label{P3.8}
The space $S(T^2)$ fails property upper $(H)$
(even for disjointly supported elements)  and fails property
lower $(H_{2})$. Hence $S(T^2)$ is not of weak type 2 and not 
of cotype 2.
\end{Prop}

\bproof
First we check property lower $(H_{2})$.  Since $(t_{n}^{s})$
is symmetric and
is dominated by the unit vector basis of ${\ell}_{2}$, it follows that
if this
family had subsets dominating the unit vector basis of 
${\ell}_{2}$, then $(t_{n}^{s})$ would be equivalent to the
unit vector basis of ${\ell}_{2}$ which is impossible. 

For property upper $(H)$, fix $M>1$ and choose a 
decreasing sequence of non-zero scalars
$(a_{i})_{i=1}^{n}$ whose ${\ell}_{2}$ norm is $>M$ but
$\|\sum_{i}a_{i}t_{i}\|_{T^2}= 1$.  This can be done by a modification
of the constructions of \cite[Chapter IV]{CS}.  Now let
$(x_{j})_{j=1}^{n}$ be a sequence of disjoint vectors in $S(T^2)$
which have this set of $a_{i}'s$ as coefficients.  So
$\|x_{i}\|_{S(T^2)}= 1$ for every $i=1,2,\ldots , n$.  But 
to norm  $\sum_{i}x_{i}$ in $S(T^2)$, we have to arrange all the
coefficients in decreasing order and take the norm in $T^2$.
Since these vectors are disjoint, at least half of them,
say $(x_{i})_{i\in I}$,  will
have all of their support after $t_{n/2}$.  That is, we have
n/2 vectors in $T^2$ which are disjoint and have their supports
after $t_{n/2}$.  Hence 
\begin{eqnarray*}
\|\sum_{i=1}^{n}x_{i}\|_{ST^2} & \ge & K^{-1} \|D\sum_{i=1}^{n}x_{i}
\|_{T^2} \ge  K^{-1}(\sum_{i\in I}\|x_{i}\|_{T^2}^{2} ) ^{1/2}\\
& \ge & K^{-2}(\sum_{i\in I}\|x_{i}\|_{ST^2}^{2})^{1/2} \ge
K^{-2}M(\frac{n}{2})^{1/2}.
\end{eqnarray*}

Since $M$ was arbitrarily large, it follows that $S(T^2)$ fails
upper $(H)$ - for disjoint elements.      
\eproof

We shall now need a result essentially due to S. Kwapien. In the form we
present it is due to W.B. Johnson and it appeared in \cite{LT}

\begin{Prop}\label{P1.1}
There is a function 
$$
N(k,{\epsilon})= \left [ \frac{2k^{2}}{\epsilon}\right ] ^{k}
$$
such that for any fixed $0<\epsilon <1$, every 
order complete Banach Lattice
L, and every k-dimensional subspace $F$ of $L$, there are $N=
N(k,{\epsilon})$ disjoint elements $(x_{j})_{j=1}^{N}$ in $L$
and a linear operator $V:F\rightarrow X=\text{span}(x_{j})$
such that for all $x\in X$ we have
$$
\|Vx-x\|\le {\epsilon}\|x\|.
$$
\end{Prop}

\begin{Prop}\label{P3.9}
There is a constant $K>1$ so that for every subspace $E$
of $S(T^2)$ of dimension n, we have for all $i\in \mathbb N$
for which  ${\text log}_{i-1}n$ exists,
$$
d(E,{\ell}_{2}^{n})\le K^{i}{\text log}_{i-2}n.
$$
\end{Prop}

\bproof
By giving up one level of logs we may assume by Proposition
\ref{P1.1} that we are working with a normalized disjointly
supported sequence of vectors $(x_{j})_{j=1}^{n}$ in $S(T^2)$.
Now there is a disjoint set of permutations $y_{j}$ of the
$x_{j}$ so that
\begin{eqnarray*}
\|\sum_{j=1}^{n}a_{j}x_{j}\|_{ST^2} & \ge &
\frac{1}{K} \|\sum_{j=1}^{n}a_{j}y_{j}\|_{T^2}\\
\ge \frac{1}{K}\|\sum_{j=1}^{n}a_{j}t_{j}\|_{T^2} & \ge &
\frac{1}{K^{i+1}
({\text log}_{i}n)}\left ( \sum_{j=1}^{n}|a_{j}|^{2}\right ) ^{1/2}.
\end{eqnarray*}
Also, let $Dx_{j}=z_{j}$ and
$$
w_{j} = \sum_{k}z_{j}(k)t_{n(k-1)+j},
$$
By Proposition \ref{P2.12} we have
\begin{eqnarray*}
\|D\sum_{j=1}^{n}a_{j}x_{j}\|_{ST^2} & \le & K\|\sum_{j=1}^{n}
a_{j}w_{j}\|_{T^2} 
\le 2K\left ( \sum_{j=1}^{n}|a_{j}|^{2}\|w_{j}\|_{T^2}^{2}
\right ) ^{1/2} \\
& \le & 2K\left ( \sum_{j=1}^{n}|a_{j}|^{2}
[2K^{i}({\text log}_{i}n)]^{2}\right ) ^{1/2} 
\le 4K^{i+1}({\text log}_{i}n)\left ( \sum_{j=1}^{n}|a_{j}|^{2}
\right ) ^{1/2},
\end{eqnarray*}
and hence

$$
d(E,{\ell}_{2}^{n})\le 4K^{2(i+1)}({\text log}_{i}n)^{2}
\le K^{i}({\text log}_{i-1}n).
$$
The ${\text log}_{i-2}n$ in the statement of the theorem
comes from the fact that we first applied Proposition \ref{P1.1}. 
\eproof

\begin{Cor}
The space $S(T^{2})$ is of type p for all $1\le p<2$ and of cotype
q for all $2<q$.  
\end{Cor}

Before we go on, we need a criterion for a Banach space to be of weak
cotype 2. We shall say that a Banach space $X$ has {\em property $(P)$}
if there is a constant K so that whenever
$\{x_1,x_2,\dots,x_n\}\subseteq X$ is a finite set with $\max_{1\le
j \le n}|t_j| \le \| \sum_{j=1}^n t_jx_j\|$ for all $(t_j) \subseteq
\bR$,
then 
\begin{equation}
\label{eq3.9}
\sqrt{n} \le K \big(\int \big\|\sum_{j=1}^n g_j(t)x_j
\big\|^2d\mu(t)\big)^\frac12
\end{equation}

It was proved by Pisier \cite[Proposition 10.8]{P} that if $X$ is of
weak
cotype 2, then it has property $(P)$. It turns out that $(P)$
characterizes weak cotype 2 spaces. This fact might be known to
specialists but we shall give a short proof here:

\begin{Thm}
\label{thm3.9a}
If $X$ has property $(P)$, then it is of weak cotype 2.
\end{Thm}

\bproof
Let $E\subseteq X$ be a finite dimensional subspace, say $\dim(E) =
2n$. By a result of Bourgain and Szarek \cite[Theorem 2]{BS} there is a
universal constant $C$ and $\{x_1, x_2, \dots,x_n\} \subseteq X$  so
that
for all $(t_j)\subseteq \bR$ we have
\begin{equation}
\label{eq3.10}
\max_{1\le j \le n} |t_j| \le \|\sum_{j=1}^n t_jx_j\| \le
C\big(\sum_{j=1}^n |t_j|^2 \big)^\frac12 
\end{equation}

Using property $(P)$ we get that
\begin{equation}
\label{eq3.11}
\sqrt{n} \le K \big(\int \big\|\sum_{j=1}^n g_j(t)x_j
\big\|^2d\mu(t)\big)^\frac12
\end{equation}
where $K$ is the constant of property $(P)$.
Now, (\ref{eq3.11}) and the right inequality of
(\ref{eq3.10}) give together with one of main results of \cite[Theorem
2.6]{FLM} (see also \cite[pages 25 and 81]{T}) that there is a universal
constant $\eta$ such that if $k\le \eta K^{-2}C^{-2} n$, then there is a
$k$-dimensinal subspace $F\subseteq [x_j]$ with $d(F,l_2^k)\le 2$. From
\cite[Theorem 10.2]{P} it now follows that $X$ is of weak cotype 2.
\eproof

We shall say that a sequence $(x_j)_{j=1}^n$ in a Banach space $X$ is
1-separated if $\|x_i - x_j\| \ge 1$ for all $1\le i,j\le n$, $i\ne
j$. It follows immediately from Theorem \ref{thm3.9a} that if every
1-separated sequence in $X$ satisfies (\ref{eq3.9}), then $X$ is of weak
cotype 2.

We are now ready to prove that the symmetric convexified Tsirelson
space is a weak cotype 2 space with a symmetric basis which is
not of cotype 2.  Hence its dual space is a symmetric space
which is of weak type 2 but fails to be of type 2.

\begin{Thm}\label{P3.10}
The space $S(T^2)$ is a weak cotype 2 space.
\end{Thm}

\bproof
Let $(x_{j})_{j=1}^{n}$ be a 1-separated
sequence in $S(T^2)$.  Without loss of generality we may assume that
for all $1\le i\le n$ we have $\|x_{i}\|_{S(T^{2})}\ge 1$. We wish to
show that (\ref{eq3.9}) holds. If $K$ is a constant which satisfies
(\ref{eq1.5}) for both $T^2$ and $S(T^2)$ and (\ref{eq3.3}), then by
definition we can find a $\sigma \in \Pi$ so that: 
\begin{equation}
\label{eq3.8}
\|(\sum_{j=1}^n |S_{\sigma}x_j|^2)^\frac12 \|_{T^2} =
\|S_{\sigma}(\sum_{j=1}^n |x_j|^2)\|_{T^2} \le K\|(\sum_{j=1}^n
|x_j|^2)^\frac12 \|_{S(T^2)}
\end{equation}

Since $S_{\sigma}$ is an isometry on $S(T^2)$, we can without loss of
generality assume that actually $x_j = S_{\sigma}x_j$ for all $ 1\le j
\le n$.

Put $k=\log \log n$ and let $P_k$ be the natural projection of $T^2$
onto
the span of $(t_j)_{j=1}^k$. We now examine two cases.

{\bf Case I}:  There is a subset $I\subset \{1,2,\cdots ,n\}$ with
$|I|\ge \frac{n}{2}$ so that $\|P_{k}x_{j}\|_{{\ell}_{2}}\ge 
\log k$ for all $j\in I$.

Since $(t_j)_{j=1}^k$ is $K\log k$-isomorphic to a Hilbert space by
Proposition \ref{P2.10}, we get using (\ref{eq1.5}) and (\ref{eq3.8})   
\begin{eqnarray}
\label{eq3.9a}
\big(\int \big\|\sum_{j=1}^n g_j(t)x_j
\big\|_{S(T^2)}^2d\mu(t)\big)^\frac12 & \ge &
\frac1K \|(\sum_{j=1}^{n}|x_{j}|^{2})^{1/2}\|_{S(T^2)}
\ge \frac{1}{K^2}\|(\sum_{j=1}^{n}|x_{j}|^{2})^{1/2}\|_{T^2}\\
\nonumber
& \ge & \frac{1}{K^2}\|(\sum_{j=1}^{n}|P_k x_{j}|^{2})^{1/2}\|_{T^2} 
\ge \frac{1}{(\log k)K^3}\|(\sum_{j\in I}|P_k x_{j}|^{2})^{1/2}
\|_{{\ell}^2} \\
\nonumber
& = & \frac{1}{(\log k)K^3} (\sum_{j\in I} \|P_k x_j
\|^2)^\frac12 \ge \frac{1}{K^3\sqrt{2}} \sqrt{n}
\end{eqnarray}

{\bf Case II}:  There is a subset $I\subset \{1,2,\cdots ,n\}$ with
$|I|\ge \frac{n}{2}$ so that $\|P_{k}x_{j}\|_{{\ell}_{2}}\le 
\log k$ for all $j\in I$.
   
In this case we make the following claim:

{\em Claim}:  There is a subset $J\subset I$
with $|J|\ge \frac{n}{4}$, so that for
all $j\in J$,
$$
\|(I-P_{k})x_{j}\|_{T^2}\ge  \frac{1}{8K}.
$$

If not, there is a set $J$ as above with
$$
\|(I-P_{k})x_{j}\|_{T^2}\le  \frac{1}{8K}.
$$

By a volume of the ball argument (see e.g. \cite[Lemma 2.4]{FLM}) the
cardinality of a set of points which are $\frac{1}{4K}$ apart in a ball
of radius $\log k$ in $k$-dimensional Hilbert space is at most
$(1+8K\log k)^k$ which by our choice of $k$ is less than or equal to
$\frac{n}{4}$ (at least for large $n$). Hence there exist $i,j\in J$,
$i \ne j$ so that
$$
\|P_{k}(x_{i}-x_{j})\|_{{\ell}^2}\le \frac{1}{4K}.
$$
Now we compute
\begin{eqnarray*}
\|x_{i}-x_{j}\|_{S(T^{2})}& \le & K\|x_{i}-x_{j}\|_{T^{2}} \le
K\|P_{k}(x_{i}-x_{j})\|_{T^{2}}+
K\|(I-P_{k})x_{i}\|_{T^{2}}+K\|(I-P_{k})x_{j}\|_{T^{2}}\\
& \le & K\|P_{k}(x_{i}-x_{j})\|_{{\ell}_{2}} + K\frac{1}{8K}
+ K\frac{1}{8K}
\le K\frac{1}{4K}+\frac{1}{4} = \frac{1}{2}.
\end{eqnarray*}
This contradicts our 1-separation assumption.  So the claim
holds.   

Now by the claim, the beginning of the proof, (\ref{eq1.5}) and
Proposition \ref{P2.11} 
we get
\begin{eqnarray}
\label{3.10a}
\big(\int \big\|\sum_{j=1}^n g_j(t)x_j
\big\|_{S(T^2)}^2d\mu(t)\big)^\frac12 & \ge &
\frac{1}{K^2}\|(\sum_{j=1}^{n}|x_{j}|^{2})^{1/2}\|_{T^2}\\
\nonumber
& \ge & \frac{1}{K^2}\|(I-P_{k})(\sum_{j\in J}|x_{j}|^{2})^{1/2}
\|_{T^2}\\ 
\nonumber
& \ge & \frac{1}{K^3}\big(\int \big\|\sum_{j=1}^n g_j(t)(I-P_k)x_j
\big\|_{T^2}^2d\mu(t)\big)^\frac12 \\
\nonumber
& \ge & \ge \frac{1}{K^5}\left ( \sum_{j\in J}\|(I-P_{k})
x_{j}\|_{T^2}^{2}\right ) ^{1/2} \\
\nonumber
& \ge & \frac{1}{K^5}\left ( \sum_{j\in J}
(\frac{1}{8K})^{2}\right ) ^{1/2} \ge \frac{|J|^{1/2}}{8 K^6} \ge
\frac{\sqrt{n}}{16K^6} 
\end{eqnarray}
This completes the proof.
\eproof

As a corollary we obtain:

\begin{Cor}
Even for Banach lattices property upper H and the weak type 2
property do not imply the upper $H_{2}$ property.  Similarly,
property lower H and the weak cotype 2 property do not imply
the lower $H_{2}$ property.
\end{Cor}

\vspace{1cm}

\noindent Department of Mathematics,\\ University of Missouri,\\
Columbia MO 65211,\\
pete@casazza.math.missouri.edu\\

\noindent Department of Mathematics and Computer Science,\\ SDU-Odense 
University,\\ Campusvej 55, DK-5230 Odense M, Denmark,\\
njn@imada.sdu.dk
  
\end{document}